\documentclass[12pt]{amsart}
\usepackage{amssymb}
\usepackage[margin=1in]{geometry}
\usepackage{hyperref}

\def\Re{\textup{Re}\,}
\def\Im{\textup{Im}\,}

\newtheorem{theorem}{Theorem}
\newtheorem{lemma}{Lemma}
\newtheorem{satz}{Satz}
\numberwithin{satz}{section}

\begin{document}

\parskip5pt
\parindent20pt
\baselineskip15pt

\title[Multiplicative Function Mean Values:  Asymptotic Estimates]{Multiplicative Function Mean Values:  \\ Asymptotic Estimates}

\author{P. D. T. A. Elliott}
\address{Department of Mathematics, University of Colorado Boulder, Boulder, Colorado 80309-0395 USA}
\email{pdtae@euclid.colorado.edu}

% THIS CODE PUTS SECTION NAMES AND THE REFERENCES ON THE LEFT:

\makeatletter
\def\specialsection{\@startsection{section}{1}%
  \z@{\linespacing\@plus\linespacing}{.5\linespacing}%
%  {\normalfont\centering}}% DELETED
  {\normalfont}}% NEW
\def\section{\@startsection{section}{1}%
  \z@{.7\linespacing\@plus\linespacing}{.5\linespacing}%
%  {\normalfont\scshape\centering}}% DELETED
  {\normalfont\scshape}}% NEW
\makeatother

% THE FOLLOWING LINE IS THE DEDICATION.  TO TURN IT ON DELETE THE % SIGN AT THE BEGINNING OF THE LINE:

%\dedicatory{In celebration of the eighty fifth birthday of Eduard Wirsing.}

%  IF NECESSARY, HERE ARE SUBJECT CLASS, KEYWORDS AND A SUMMARY:

%\subjclass[2010]{Primary 11N37; Secondary 11F03, 11F30, 11K65, 11L99, 11M99, 11N60, 11N64}

%\keywords{Multiplicative function, Mean value, Maass form}

%\begin{abstract}
%
%Mean value theorems for multiplicative arithmetic functions are applied to demonstrate uniformity of sign changes in the Fourier coefficients of automorphic forms.
%
%\end{abstract}

\maketitle

\section{Statement of results}

For many studies in analytic number theory a natural object against which to measure the mean-value of a complex-valued multiplicative arithmetic function $n \to g(n)$ is the mean-value of its attendant function $n \to |g(n)|$.

This reflects the decomposition $n \to |g(n)| \exp( i \arg g(n))$ of a non-vanishing completely multiplicative function into essentially a unitary character on the multiplicative group of the positive rationals, and a homomorphism $n \to \log |g(n)|$ of the positive rationals into the additive reals.

Some fifty years ago, papers of Delange \cite{delange1961surlesfonctionsapplications} 1961, Wirsing \cite{wirsing1961} 1961, \cite{wirsing1967} 1967, Hal\'asz \cite{halasz1968} 1968, catalysed the general study of multiplicative functions and moved the field seriously forward.

In the present paper I re-examine the theorems of Wirsing in the light of more recent developments and apply related ideas to the consideration of two open-ended questions.

The following four cumulative theorems will be established, all new.  Several auxiliary propositions are also of independent interest.

\begin{theorem} \label{elliott_2015_multiplicative_function_mean_values_thm_01}
Let $g$ be a non-negative multiplicative function, uniformly bounded on the primes, for which the series $\sum q^{-1} g(q)$, taken over the prime-powers $q = p^k$ with $k \ge 2$, converges, and for which the sums $y^{-1} \sum_{q \le y} g(q)\log q$, $y \ge 2$, are uniformly bounded.

Let $h(n)$ be a complex-valued multiplicative function that satisfies $|h(n)| \le g(n)$.

Set $G(x) = \sum_{n \le x} n^{-1} g(n)$, $H(x) = \sum_{n \le x} n^{-1} h(n)$, $x \ge 1$.

Then
\[
H(x) = \left( \prod\limits_{p \le x} \left( 1 + \frac{h(p)}{p} + \frac{h(p^2)}{p^2} + \cdots \right) \left( 1 + \frac{g(p)}{p} + \frac{g(p^2)}{p^2} + \cdots \right)^{-1} + o(1) \right) G(x)
\]
as $x \to \infty$.
\end{theorem}

\noindent \emph{Remark}.  If the series $\sum p^{-1} ( g(p) - \Re h(p))$ diverges or a sum $\sum_{k=1}^\infty p^{-k} h(p^k)$ has the value $-1$, then the product over the primes may be omitted.  Otherwise, the product has the form $AL(\log x)$, where $A$ is a non-zero constant and $L(y)$ a non-vanishing slowly oscillating function of $y$.

\begin{theorem} \label{elliott_2015_multiplicative_function_mean_values_thm_02}
Let $g$ be a non-negative multiplicative function that is uniformly bounded on the primes.  Assume that for a positive $c$, and each $b$, $0 < b < 1$,
\[
\liminf \limits_{x \to \infty} ( (1-b)\log x)^{-1} \sum\limits_{x^b < p \le x} p^{-1} g(p)\log p \ge c.
\]

Then for some positive $c_0$ and all $x \ge 2$,
\[
\sum\limits_{n \le x} g(n) \ge \frac{c_0x}{\log x} \prod\limits_{p \le x} \left( 1 + \frac{g(p)}{p} \right).
\]
\end{theorem}

\noindent \emph{Remark}.  Under the further assumptions on $g$ in Theorem \ref{elliott_2015_multiplicative_function_mean_values_thm_01}, there is a similar upper bound.

For each positive real $\tau$, $\Delta(\tau)$ will denote a compact star-shaped region of the complex plane that contains the origin, has a representation
\[
\{ \rho e^{i\theta}, 0 \le \theta < 2\pi, 0 \le \rho \le w(\theta)\},
\]
with average radius
\[
(2\pi)^{-1} \int_0^{2\pi} w(\theta) \, d\theta, \quad w(2\pi) = w(0),
\]
strictly less than $\tau$.

\begin{theorem} \label{elliott_2015_multiplicative_function_mean_values_thm_03}
Let the multiplicative function $g$ satisfy the hypotheses of Theorems \ref{elliott_2015_multiplicative_function_mean_values_thm_01} and \ref{elliott_2015_multiplicative_function_mean_values_thm_02} and let $h$ be a complex-valued multiplicative function with $|h(n)| \le g(n)$ and values in $\Delta(c)$.

Set
\[
A(x) = \sum\limits_{n \le x} g(n), \quad B(x) = \sum\limits_{n \le x} h(n).
\]

Then
\[
B(x) =  \left( \prod\limits_{p \le x} \left( 1 + \frac{h(p)}{p} + \frac{h(p^2)}{p^2} + \cdots \right) \left( 1 + \frac{g(p)}{p} + \frac{g(p^2)}{p^2} + \cdots \right)^{-1} + o(1) \right) A(x)
\]
as $x \to \infty$.
\end{theorem}

\begin{theorem} \label{elliott_2015_multiplicative_function_mean_values_thm_04}
Let the multiplicative function $g$ satisfy the hypotheses of Theorems \ref{elliott_2015_multiplicative_function_mean_values_thm_01} and \ref{elliott_2015_multiplicative_function_mean_values_thm_02} and let $h$ be a complex-valued multiplicative function with $|h(n)| \le g(n)$.

Then there are two possibilities.

\begin{enumerate}
\item[(i)]  For some real $t$ the series $\sum p^{-1} ( g(p) - \Re h(p)p^{it})$, taken over the primes, converges;
\[
B(x) = (1-it)^{-1} x^{-it} \prod\limits_{p \le x} \left( 1 + h(p)p^{it-1} + \cdots \right) \left( 1 + g(p)p^{-1} + \cdots \right)^{-1} A(x) + o(A(x)), \quad x \to \infty.
\]
\item[(ii)]  There is no such $t$, and
\[
B(x) = o(A(x)), \quad x \to \infty.
\]
\end{enumerate}
\end{theorem}

Of particular interest in Theorems \ref{elliott_2015_multiplicative_function_mean_values_thm_01} and \ref{elliott_2015_multiplicative_function_mean_values_thm_04} is that beyond dominance by $g$, there is no non-structural constraint upon the complex values of the function $h$.

%@@@@@@@@@@@@@@@@@@@@@@@@@@@@@@@@@@@@@@@@@@@@@@@@@@@@@@@@@

\section{Background} \label{elliott_2015_multiplicative_function_mean_values_sec_02}

Two central theorems of Wirsing's 1967 paper run as follows.

\setcounter{section}{1}

\begin{satz} \label{wirsing_1967_satz_01_01}
Let $\lambda(n)$ be a non-negative multiplicative function, uniformly bounded on the primes, that for a positive $\tau$ satisfies
\[
\sum\limits_{p \le x} p^{-1}\log p\, \lambda(p) \sim \tau \log x, \quad x \to \infty. 
\]
Assume further that the series $\sum q^{-1} \lambda(q)$, taken over the prime-powers $q = p^k$ with $k \ge 2$, converges, and that if $\tau \le 1$ then $\sum_{q \le x} \lambda(q) \ll x(\log x)^{-1}$ holds for $x \ge 2$.

Then
\[
\sum\limits_{n \le x} \lambda(n) \sim \frac{e^{-\gamma \tau}}{\Gamma(\tau)} \frac{x}{\log x} \prod\limits_{p \le x} \left(1 + \frac{\lambda(p)}{p} + \frac{\lambda(p^2)}{p^2} + \cdots \right), \quad x \to \infty, 
\]
where $\gamma$ is Euler's constant.
\end{satz}

%@@@@@@@@@@@@@@@@@@@@@@@@@@@@@@@@@@@@@@@@@@@@@@@@@@@@@@@@@

\begin{satz} \label{wirsing_1967_satz_01_02_02}
Let $\lambda(n)$ be a multiplicative function that satisfies the conditions of Satz \ref{wirsing_1967_satz_01_01}.  Let $\lambda^*(n)$ be multiplicative, with values in $\Delta(\tau)$ and satisfy $|\lambda^*(n)| \le \lambda(n)$.

Then
\[
 \sum\limits_{n \le x} \lambda^*(n) = \frac{e^{-\gamma \tau}}{\Gamma(\tau)} \frac{x}{\log x} \prod\limits_{p\le x} \left( 1 + \frac{\lambda^*(p)}{p} + \frac{\lambda^*(p^2)}{p^2} + \cdots \right) + o\left( \sum\limits_{n \le x} \lambda(n) \right)
\]
as $x \to \infty$.
\end{satz}

\setcounter{section}{2}

In what follows, a product of the form
\[
\prod\limits_{p \le x} \left( 1 + \frac{f(p)}{p} + \frac{f(p^2)}{p^2} + \cdots \right),
\]
when meaningful, may be denoted by $\sideset{}{_x}\prod  (f)$.

The two theorems of Wirsing may be compared to the following result of Elliott and Kish \cite{elliottkish2013harmonicmaass}, subsuming ideas from Wirsing and Hal\'asz, loc, cit.

%@@@@@@@@@@@@@@@@@@@@@@@@@@@@@@@@@@@@@@@@@@@@@@@@@@@@@@@@@

\begin{theorem} \label{elliott_2015_multiplicative_function_mean_values_thm_05}
Let $3/2 \le Y \le x$.  Let $g$ be a complex-valued multiplicative function that for positive constants $\beta$, $c$, $c_1$ satisfies $|g(p)| \le \beta$,
\[
\sum\limits_{w < p \le x} p^{-1} ( |g(p)| - c) \ge -c_1, \quad Y \le w \le x,
\]
on the primes.  Suppose, further, that the series
\[
\sum\limits_{q} q^{-1} |g(q)| (\log q)^\kappa, \quad \kappa = 1 + c\beta(c+\beta)^{-1},
\]
taken over the prime-powers $q = p^k$ with $k \ge 2$, converges.

Then with
\[
\lambda = \min\limits_{|t| \le T} \sum\limits_{Y < p \le x} p^{-1} ( |g(p)| - \Re g(p) p^{it}),
\]
\[
\sum\limits_{n \le x} g(n) \ll x(\log x)^{-1} \prod\limits_{p \le x} ( 1 + p^{-1} |g(p)| ) \left( \exp( -\lambda c(c+\beta)^{-1}) + T^{-1/2} \right)
\]
uniformly for $Y, x, T>0$, the implied constant depending at most upon $\beta$, $c$, $c_1$ and a bound for the sum of the series over higher prime-powers.
\end{theorem}

An extension of Theorem \ref{elliott_2015_multiplicative_function_mean_values_thm_05}, a proof of which will be given following that for Theorem \ref{elliott_2015_multiplicative_function_mean_values_thm_04}, obviates the awkward condition involving the factor $(\log q)^\kappa$.

%@@@@@@@@@@@@@@@@@@@@@@@@@@@@@@@@@@@@@@@@@@@@@@@@@@@@@@@@@

\begin{theorem} \label{elliott_2015_multiplicative_function_mean_values_thm_06}
If the estimate in Theorem \ref{elliott_2015_multiplicative_function_mean_values_thm_05} is weakened to
\[
\sum\limits_{n \le x} g(n) \ll x(\log x)^{-1} \sideset{}{_x}\prod  (|g|) \left( \exp( -\lambda c(c+\beta)^{-1}) + T^{-1/2} \right)^{c/(3c+1)},
\]
then the condition on the prime-power values $g(p^k)$, $k \ge 2$, may be relaxed to the convergence of the series $\sum_{p, k \ge 2} p^{-k} |g(p^k)|$ and a uniform bound for the sums $y^{-1} \sum_{p^k \le y} |g(p^k)| \log p^k$, $y \ge 2$.
\end{theorem}

For the multiplicative function $\lambda_0(n)$ defined to be $\alpha$, $\beta$ with $0 < \alpha < \beta$, on the primes in alternate intervals $(\exp(2^k), \exp(2^{k+1})]$, $k = -1, 0, 1, 2, \dots,$ and to be zero on all other prime-powers,
\[
\lim\limits_{x \to \infty} (\log x)^{-1} \sum\limits_{p \le x} p^{-1} \lambda_0(p) \log p
\]
does not exist, eliminating direct application of S\"atze \ref{wirsing_1967_satz_01_01} and \ref{wirsing_1967_satz_01_02_02}.

The lower bound of Theorem \ref{elliott_2015_multiplicative_function_mean_values_thm_02} is obtained in Elliott and Kish \cite{elliottkish2013harmonicmaass}, Lemma 21, subject to the existence of a positive constant $c_2$ so that for all large $x$, $\sum_{p \le x} g(p) \log p \ge c_2 x$.  By modifying $\lambda_0$ to be zero on intervals $(y(\log y)^{-2}, y]$, $y = \exp(2^k)$, we obtain a multiplicative function $\lambda_1$ that will not satisfy such a criterion for any positive $c_2$.

Never-the-less, Theorems \ref{elliott_2015_multiplicative_function_mean_values_thm_01}, \ref{elliott_2015_multiplicative_function_mean_values_thm_02} and \ref{elliott_2015_multiplicative_function_mean_values_thm_04} may be applied to $\lambda_0$, $\lambda_1$ with any dominated complex-valued multiplicative function, $h$.

%@@@@@@@@@@@@@@@@@@@@@@@@@@@@@@@@@@@@@@@@@@@@@@@@@@@@@@@@@

\section{Proof of Theorem \ref{elliott_2015_multiplicative_function_mean_values_thm_01}} 

It is convenient to introduce several preliminary results.

\begin{lemma} \label{elliott_2015_multiplicative_function_mean_values_lem_01}
The estimate
\[
\sum\limits_{2 \le n \le x} g(n) \le \left( \frac{x}{\log x} + \frac{10x}{(\log x)^2} \right) \widetilde{\Delta} \sum\limits_{n \le x} \frac{g(n)}{n}
\]
with
\[
\widetilde{\Delta} = \sup\limits_{1 \le y \le x} y^{-1} \sum\limits_{q \le y} g(q)\log q,
\]
where $q$ denotes a prime-power, holds uniformly for all non-negative real multiplicative functions $g$, and all $x \ge 2$.
\end{lemma}

A proof of Lemma \ref{elliott_2015_multiplicative_function_mean_values_lem_01} may be found in Elliott \cite{Elliott1997}, Chapter 2, Lemma 2.2.  It is immediate that
\begin{align*}
\sum\limits_{n \le x} n^{-1} g(n) & \le \prod\limits_{p \le x} \left(1 + \sum\limits_{k \le \log x/\log p} p^{-k} g(p^k) \right) \\
& \le \exp\left( \sum\limits_{q \le x} q^{-1} g(q) \right).
\end{align*}

A proof of the following qualitative corresponding lower bound,  a result first obtained by Barban \cite{barban1966largesieve} using a different method, may be found in Lemma 20 of Elliott and Kish, \cite{elliottkish2013harmonicmaass}.

%@@@@@@@@@@@@@@@@@@@@@@@@@@@@@@@@@@@@@@@@@@@@@@@@@@@@@@@@@

\begin{lemma} \label{elliott_2015_multiplicative_function_mean_values_lem_02}
To each positive $\beta$ there is a further positive $c(\beta)$ so that a non-trivial non-negative multiplicative function, $g$, that satisfies $g(p) \le \beta$ on the primes, also satisfies 
\[
\sum\limits_{n \le x} g(n) n^{-1} \ge c(\beta) \prod\limits_{p \le x} ( 1 + p^{-1} g(p) )
\]
uniformly for $x \ge 1$.
\end{lemma}

%@@@@@@@@@@@@@@@@@@@@@@@@@@@@@@@@@@@@@@@@@@@@@@@@@@@@@@@@@

\begin{lemma} \label{elliott_2015_multiplicative_function_mean_values_lem_03}
Let $g$ be a non-trivial non-negative multiplicative function uniformly bounded on the primes, for which the series $\sum q^{-1} g(q)$, taken over the prime-powers $q = p^k$ with $k \ge 2$, converges, and for which the sums $y^{-1} \sum_{q \le y} g(q) \log q$, $y \ge 2$, are uniformly bounded.

Then
\[
\sum\limits_{u < n \le v} \frac{g(n)}{n} \ll \left( \log \left( \frac{\log v}{\log u}\right) + \frac{1}{\log x} \right) \sum\limits_{n \le x} \frac{g(n)}{n}
\]
uniformly for $x^{1/2} \le u \le v \le x^{3/2}$, $x \ge 2$.
\end{lemma}

%@@@@@@@@@@@@@@@@@@@@@@@@@@@@@@@@@@@@@@@@@@@@@@@@@@@@@@@@@

\noindent \emph{Proof of Lemma \ref{elliott_2015_multiplicative_function_mean_values_lem_03}}.  In view of the hypothesis on $g$, Lemma \ref{elliott_2015_multiplicative_function_mean_values_lem_01} delivers the uniform estimate
\[
\sum\limits_{n \le y} g(n) \ll \frac{y}{\log y} \prod\limits_{p \le x^{3/2}} \left( 1 + \frac{g(p)}{p} \right), \quad 2 \le y \le x^{3/2},
\]
which Lemma \ref{elliott_2015_multiplicative_function_mean_values_lem_02} shows to be $\ll y(\log y)^{-1} G(x)$.  The asserted result then follows from an integration by parts.

For better appreciation the following theorem is given in both its abelian and tauberian aspects.  A proof may be found, together with a history of the result from Feller \cite{feller1963tauberian} to Stadtm\"uller and Trautner \cite{stadtmullertrautner1979tauberianlaplace}, in Bingham, Goldie and Teugels \cite{binghametal1987regularvariation}, Chapter 2, Theorem 2.10.1, pp. 116--118, and Korevaar \cite{korevaar2004tauberiantheory}, Chapter IV, Theorem 10.1, pp. 197--199.

Let $C(y)$, $D(y)$ be non-negative real-valued functions on the non-negative reals, non-decreasing and right continuous.  To each corresponds a Laplace transform, typically
\[
s \to \widehat{C}(s) = \int_0^\infty e^{-sy} \, dC(y),
\]
here assumed to be defined for $s > 0$.

%@@@@@@@@@@@@@@@@@@@@@@@@@@@@@@@@@@@@@@@@@@@@@@@@@@@@@@@@@

\begin{lemma} \label{elliott_2015_multiplicative_function_mean_values_lem_04}
Assume that for each $y > 1$
\[
D^*(y) = \limsup \limits_{u \to \infty} D(u)^{-1} D(uy)
\]
is finite, $D$ implicitly assumed not to be identically zero.

If, for some constant $A$ and slowly-oscillating function $L(y)$,
\[
C(y) = (AL(y) + o(1)) D(y), \quad y \to \infty,
\]
then
\[
\widehat{C}(s) = (AL(s^{-1}) + o(1)) \widehat{D}(s), \quad s \to 0+.
\]

Further, if $D^*(y) \to 1$ as $y \to 1+$, then the converse is valid.
\end{lemma}

\noindent \emph{Remark}. The non-decreasing nature of $D$ ensures that $\lim D^*(y)$, $y \to 1$, exists.

\noindent \emph{Completion of the proof of Theorem \ref{elliott_2015_multiplicative_function_mean_values_thm_01}}.  We apply Lemma \ref{elliott_2015_multiplicative_function_mean_values_lem_04} to the pair $2G(e^x) + \Re( H(e^x))$, $G(e^x)$; to the pair with $\Im(H(e^x))$ in place of $\Re(H(e^x))$; and to the pair $G(e^x)$, $G(e^x)$.

Computation with Euler products shows $\widehat{C}(s)$, $\widehat{D}(s)$, the Laplace transforms of the first pair, to exist for all positive $s$ and satisfy $\widehat{C}(s) = f(s)\widehat{D}(s)$, where
\[
f(s) - 2 = \Re \left( \prod\limits_{p} \left( 1 + \sum\limits_{k=1}^\infty p^{-k(1+s)}h(p^k) \right)\left( 1 + \sum\limits_{m=1}^\infty p^{-m(1+s)} g(p^m) \right)^{-1} \right).
\]
In particular,
\[
|f(s) - 2| \ll \exp\left( - \sum_p p^{-1-s}( g(p) - \Re h(p)) \right),
\]
so that if the series in the exponent diverges for $s=0$, then $f(s) \to 2$ as $s \to 0+$, and we may apply Lemma \ref{elliott_2015_multiplicative_function_mean_values_lem_04} with $A=2$, $L$ identically 1.

We may therefore assume the series $\sum p^{-1}( g(p) - \Re h(p))$ to converge.

From the Chebyshev bound $\pi(y) \ll y(\log y)^{-1}$, integration by parts shows the series $\sum_{p > x^\varepsilon} p^{-1} \exp( -\log p / \log x)$ to be bounded in terms of $\varepsilon$ alone.  Since
\[
|g(p) - h(p)|^2 \le 2g(p) (g(p) - \Re h(p)),
\]
an application of the Cauchy-Schwarz inequality, confined to the primes on which $g$ does not vanish, shows that
\[
\sum\limits_{p > x^\varepsilon} p^{-1 - 1/\log x} |g(p) - h(p)| \ll \left( \sum\limits_{p > x^\varepsilon} g(p) p^{-1-1/\log x} \right)^{1/2} \left( \sum\limits_{p > x^\varepsilon} p^{-1} (g(p) - \Re h(p)) \right)^{1/2}
\]
and $o(1)$ as $x \to \infty$.

Moreover,
\[
\sum\limits_{p \le x^\varepsilon} ( p^{-1} - p^{-1 - 1/\log x}) \ll \sum\limits_{p \le x^\varepsilon} p^{-1} \log p / \log x \ll \varepsilon,
\]
the implied constant absolute for all values of $x$ sufficiently large in terms of $\varepsilon$.

Letting $x \to \infty$, $\varepsilon \to 0+$, we see that as $x \to \infty$
\[
f\left( (\log x)^{-1} \right) - 2 = \Re \left( B \exp \left( \sum\limits_{p \le x} p^{-1} \Im ( h(p)) \right)\right) + o(1),
\]
with $B$ the product of 
\[
\prod_p \left( 1 + \sum_{k=1}^\infty p^{-k} h(p^k) \right) \exp\left( -p^{-1} h(p) \right) \prod_p \left( 1 + \sum_{m=1}^\infty p^{-m} g(p^m) \right)^{-1} \exp \left( p^{-1}g(p) \right)
\]
and $\exp(-\sum_p p^{-1} (g(p) - \Re h(p)) )$.  Its genesis in terms of Euler products ensures that $|B| \le 1$; moreover, $B$ will vanish only if for some prime $p$ the sum $1 + \sum_{k=1}^\infty p^{-k} h(p^k)$ vanishes.

Note that for any $\beta \ge 1$, the above argument shows that
\[
\sum_{x < p \le x^\beta} p^{-1} \Im (h(p)) = - \sum_{x < p \le x^\beta} p^{-1} \Im (g(p) - h(p))
\]
\[
\ll \left( \sum_{x < p \le x^\beta} p^{-1} \right)^{1/2} \left( \sum_{p > x} p^{-1} |g(p) - h(p)|^2 \right)^{1/2} = o(1)
\]
as $x \to \infty$, so that $\exp( \sum_{p \le e^s} p^{-1} \Im( h(p)))$ is a slowly oscillating function of $s$.

In view of Lemma \ref{elliott_2015_multiplicative_function_mean_values_lem_03},
\[
\lim\limits_{y \to 1+} \limsup\limits_{u \to \infty} G(e^x)^{-1} G(e^{xy}) = 1.
\]
Three applications of Lemma \ref{elliott_2015_multiplicative_function_mean_values_lem_04} in its Tauberian aspect, typically with $A=1$,
\[
L(s) = 2 + \Re \left( B \exp\left( \sum_{p \le e^s} p^{-1} \Im(h(p)) \right) \right),
\]
delivers the asymptotic estimate
\[
H(e^x) = \left( f(x^{-1}) + o(1)\right) G(e^x), \quad x \to \infty,
\]
from which Theorem \ref{elliott_2015_multiplicative_function_mean_values_thm_01} follows rapidly.

%@@@@@@@@@@@@@@@@@@@@@@@@@@@@@@@@@@@@@@@@@@@@@@@@@@@@@@@@@

\section{Proof of Theorem \ref{elliott_2015_multiplicative_function_mean_values_thm_02}} \label{elliott_2015_multiplicative_function_mean_values_sec_04}

Again a preliminary result is advantageous.

Let $0 \le g(p) \le \beta$ for each prime, $p$.

If, for some $\tau > 0$,
\[
\sum_{p \le y} p^{-1} g(p) \log p \sim \tau \log y, \quad y \to \infty,
\]
then for each $\varepsilon$, $0 < \varepsilon < 1$,
\[
\liminf\limits_{x \to \infty} (\varepsilon \log x)^{-1} \sum_{x^{1-\varepsilon} < p \le x} p^{-1} \log p \ge \tau.
\]

The converse need not be true, as may be seen from the example $\lambda_0$ in section \ref{elliott_2015_multiplicative_function_mean_values_sec_02}.

However, the following converse is valid.

%@@@@@@@@@@@@@@@@@@@@@@@@@@@@@@@@@@@@@@@@@@@@@@@@@@@@@@@@@

\begin{lemma} \label{elliott_2015_multiplicative_function_mean_values_lem_05}
Assume that for $c>0$ and each $\varepsilon$, $0 < \varepsilon < 1$, the function $g(p)$, uniformly bounded on the primes, satisfies
\[
\liminf \limits_{x \to \infty} (\varepsilon \log x)^{-1} \sum_{x^{1-\varepsilon} < p \le x} p^{-1} g(p)\log p \ge c.
\]
Then for each $\alpha$, $0 < \alpha < c$, there is a subsequence of primes, $r$, such that
\[
\lim\limits_{x \to \infty} (\log x)^{-1} \sum_{r \le x} r^{-1} g(r) \log r = \alpha.
\]
\end{lemma}

\noindent \emph{Proof of Lemma \ref{elliott_2015_multiplicative_function_mean_values_lem_05}}.  We begin with an outline of the argument.  Fix a prime $t$ for which $\sum_{p \le t} p^{-1} g(p) \log p \ge \alpha \log t$.

We define a function $\overline{g}(p)$ by choosing, for each prime $p$, to retain $g(p)$ or to replace it by zero.  For ease of notation $\sum_{p \le y} p^{-1} \overline{g}(p) \log p$ will be denoted by $S(y)$.

We choose $\overline{g}(p) = g(p)$ for $p \le t$.

The primes $y_1 < y_2 < \cdots$ are defined successively as follows.  We replace $g(p)$ by zero on the primes following $t$ until, for the first time, $S(y)/\log y$ falls strictly below $\alpha$.  The corresponding value of $y$ is $y_1$.

We choose $\overline{g}(p) = g(p)$ on the primes $p > y_1$ until, for the first time with $y > y_1$, the ratio $S(y)/\log y$ climbs above $\alpha$.  The corresponding value of $y$ is $y_2$; and so on.

Our initial aim is to show the turning values $y_j$ not to be logarithmically far apart.

A few preliminary remarks are helpful.

Let $0 < \theta < 1$, $x \ge 2$, $3/2 \le y \le x^\theta$.  With $0 < \varepsilon < 1-\theta$ determine the integer $k$ by $x^{(1-\varepsilon)^k} < y \le x^{(1-\varepsilon)^{k-1}} = \psi$, so that $k \ge 2$.  Assume that for all sufficiently large values of $w$
\[
\sum_{w^{1-\varepsilon} < p \le w} p^{-1} g(p) \log p \ge \varepsilon c \log w.
\]
By partitioning the interval $(x^{(1-\varepsilon)^k}, x]$ into adjoining subintervals $(x^{(1-\varepsilon)^m}, x^{(1-\varepsilon)^{m-1}}]$, $m = 1, 2, \dots, k$, we see that provided $x^{(1-\varepsilon)^k}$ is sufficiently large in terms of $\varepsilon$,
\begin{align*}
\sum_{y < p \le x} p^{-1} g(p) \log p & \ge c \log(x/\psi) \ge c( \log(x/y) - \log(\psi/y)) \\
& \ge c\left( 1 - \varepsilon(1-\theta)^{-1} \right) \log(x/y),
\end{align*}
since $\log(\psi/y) \le \log(\psi/\psi^{1-\varepsilon}) = \varepsilon \log \psi \le \varepsilon \log x \le \varepsilon (1-\theta)^{-1} \log(x/y)$.

For the purposes of proving Lemma \ref{elliott_2015_multiplicative_function_mean_values_lem_05} we may therefore replace its lower-bound hypothesis by:

\emph{For each $\varepsilon$, $0 < \varepsilon < 1$,
\[
\sum_{y < p \le x} p^{-1} g(p)\log p \ge c\log(x/y)
\]
uniformly for $1 \le y \le x^{1-\varepsilon}$ and all $x$ sufficiently large in terms of $\varepsilon$.}

It is clear that the initial prime $t$ exists.

As a second preliminary remark, if $2 \le y \le w$, then
\[
(\log w)^{-1} S(w) - (\log y)^{-1} S(y) = \left( (\log w)^{-1} - (\log y)^{-1} \right)S(y) + (\log w)^{-1} \sum_{y < p \le w} p^{-1} \overline{g}(p) \log p.
\]
Hence
\[
\left| (\log w)^{-1} S(w) - (\log y)^{-1} S(y) \right| \le (\log w \log y)^{-1} S(y)\log(w/y)
\]
\[
+ c_0(\log w)^{-1} \sum_{y < p \le w} p^{-1}\log p \le c_1 ( \log(w/y) + 1) (\log w)^{-1}
\]
with a positive constant $c_1$ dependent at most upon the upper bound for the $g(p)$.  Here we have employed the elementary estimate $\sum_{p \le y} p^{-1} \log p = \log y + O(1)$, $y \ge 2$.

In particular, if $y$ is a prime adjacent to a turning value $y_k$, then
\[
S(y)/\log y - S(y_k)/\log y_k \ll ( | \log(y_k/y)| + 1 )/ \log y_k \ll 1/\log y_k,
\]
since the ratio of successive increasing primes approaches 1.

We now show the $y_j$ not to increase too rapidly.

Suppose that $S(y_k)/\log y_k < \alpha$, so that for the next prime $p > y_k$, $g(p)$ is kept.  In particular $S(y_k) \ge \alpha \log y_k + O(1)$.  If $y_k < (\tfrac{1}{2}y_{k+1})^{1-\varepsilon} < \tfrac{1}{2} y_{k+1}$ and $y_k$ is sufficiently large, then $\tfrac{1}{2} y_{k+1} y_k^{-1} > y_k^\varepsilon$,
\begin{align*}
S(\tfrac{1}{2} y_{k+1}) & = S(\tfrac{1}{2} y_{k+1}) - S(y_k) + S(y_k) \\
& \ge c\log( \tfrac{1}{2} y_{k+1}y_k^{-1}) + \alpha \log y_k + O(1) \\
& = \alpha \log(\tfrac{1}{2}y_{k+1}) + (c-\alpha)\log( \tfrac{1}{2} y_{k+1}y_k^{-1}) + O(1).
\end{align*}
With $w$ a nearest prime to $\tfrac{1}{2}y_{k+1}$, $S(w)/\log w > \alpha$ before the next change point, $y_{k+1}$.

Thus $y_k \ge (\tfrac{1}{2}y_{k+1})^{1-\varepsilon}$.

If $S(y_k) \ge \alpha \log y_k$, then again $S(y_k) = \alpha \log y_k + O(1)$, and $\overline{g}(p) = 0$ on the primes in the interval $(y_k, \tfrac{1}{2}y_{k+1}]$.  Hence
\begin{align*}
S(\tfrac{1}{2}y_{k+1}) ( \log( \tfrac{1}{2}y_{k+1}))^{-1} & = S(y_k) (\log( \tfrac{1}{2}y_{k+1}))^{-1} \\
& = \alpha \log y_k (\log y_{k+1})^{-1} + O((\log y_k)^{-1}).
\end{align*}
If, now, $y_k < y_{k+1}^{1-\varepsilon}$ and $y_k$ is sufficiently large then
\[
S(\tfrac{1}{2}y_{k+1}) (\log (\tfrac{1}{2} y_{k+1}))^{-1} \le \alpha (1-\varepsilon) + O((\log y_k)^{-1}),
\]
again leading to a premature change point.

In this case $y_k \ge y_{k+1}^{1-\varepsilon}$.

For all large values of $y_k$, $\tfrac{1}{2}y_{k+1}^{1-\varepsilon} \le y_k \le y_{k+1}$.  As a consequence
\[
S(y_{k+1})/ \log y_{k+1} - S(y_k)/\log y_k \ll \log( y_{k+1}/y_k) / \log y_{k+1} \ll \varepsilon,
\]
the implied constant independent of $\varepsilon$.  Since $S(y_k)/\log y_k = \alpha + O(1/\log y_k)$, $S(y)/\log y - \alpha \ll \varepsilon$ for all sufficiently large values of $y$, first for prime values then for otherwise arbitrary real values.

The construction of the function $\overline{g}$ does not depend upon the value of $\varepsilon$ and we may apply the argument with $\varepsilon = 2^{-m}$, $m = 1, 2, 3, \dots$, in turn.

Lemma \ref{elliott_2015_multiplicative_function_mean_values_lem_05} is established.

%@@@@@@@@@@@@@@@@@@@@@@@@@@@@@@@@@@@@@@@@@@@@@@@@@@@@@@@@@

\emph{Completion of the proof of Theorem \ref{elliott_2015_multiplicative_function_mean_values_thm_02}}.

Let $0 < \alpha < c$ and let $r$ run through a sequence of primes for which $\sum_{r \le y} r^{-1} g(r)\log r \sim \alpha \log y$, $y \to \infty$.

Define multiplicative functions $g_j$, $j = 1, 2$, by
\[
g_1(p) = \begin{cases} g(p) & \text{if} \ p \ne r, \\ 0 & \text{if} \ p=r, \end{cases} \qquad g_2(p) = \begin{cases} 0 & \text{if} \ p \ne r, \\ g(p) & \text{if} \ p = r, \end{cases}
\]
and $g_j(p^k)=0$ on all other prime powers.

On squarefree integers $g$ coincides with $g_1 * g_2$, the Dirichlet convolution of $g_1$ and $g_2$; hence
\[
\sum_{n \le x} g(n) \ge \sum_{u \le \sqrt{x}} g_1(u) \sum_{v \le x/u} g_2(v).
\]
Satz \ref{wirsing_1967_satz_01_01} of Wirsing (c.f. \S \ref{elliott_2015_multiplicative_function_mean_values_sec_02}) gives for a typical innersum the asymptotic estimate
\[
( 1 + o(1)) \frac{e^{-\gamma \alpha}}{\Gamma(\alpha)} \frac{x}{u \log(x/u)} \prod_{p \le x/u} \left( 1 + \frac{g_2(p)}{p} \right), \quad x/u \to \infty.
\]
The doublesum thus exceeds a constant multiple of
\[
\frac{x}{\log x} \prod_{p \le x} \left( 1 + \frac{g_2(p)}{p} \right) \sum_{u \le \sqrt{x}} \frac{g_1(u)}{u}.
\]

An appeal to Lemma \ref{elliott_2015_multiplicative_function_mean_values_lem_02} completes the proof.

%@@@@@@@@@@@@@@@@@@@@@@@@@@@@@@@@@@@@@@@@@@@@@@@@@@@@@@@@@

\section{Proof of Theorem \ref{elliott_2015_multiplicative_function_mean_values_thm_03}}\label{elliott_2015_multiplicative_function_mean_values_sec_05}

Choose a real $\alpha$ to lie strictly between the average radius of $\Delta(c)$, and $c$.

Choose a subsequence of primes $r$ for which
\[
\sum_{r \le y} r^{-1} g(r)\log r \sim \alpha \log y, \quad y \to \infty.
\]

We define multiplicative functions $g_j$, $j = 1, 2$, by
\[
g_1(p^k) = \begin{cases} g(p^k) & p \ne r, \\ 0 & \text{otherwise} \end{cases}, \qquad g_2(p^k) = \begin{cases} g(p^k) & p = r, \\ 0 & \text{otherwise.} \end{cases}
\]
The function $g$ has a Dirichlet convolution representation $g_1 * g_2$.

We likewise define multiplicative functions $h_j$, $j = 1, 2$, so that $h = h_1 * h_2$, $|h_j| \le g_j$, $j= 1, 2$.  There is a representation
\[
M = \sum_{n \le x} h(n) = \sum_{u \le x} h_1(u) \sum_{v \le x/u} h_2(v).
\]

Let $0 < \varepsilon < 1/2$.  We remove the contribution from the terms with $u \le x^\varepsilon$ and $x^{1-\varepsilon} < u \le x$.  Typically, by Lemma \ref{elliott_2015_multiplicative_function_mean_values_lem_01},
\begin{align*}
\sum_{u \le x^\varepsilon} g_1(u) \sum_{v \le x/u} g_2(v) & \ll \sum_{u \le x^\varepsilon} g_1(u) \frac{x}{u \log (x/u)} \prod_{p \le x/u} \left( 1 + \frac{g_2(p)}{p} + \cdots \right) \\
& \ll \frac{x}{\log x} \sideset{}{_x}\prod  (g_2) \sum_{u \le x^\varepsilon} \frac{g_1(u)}{u}.
\end{align*}
Moreover,
\[
\sum_{u \le x^\varepsilon} \frac{g_1(u)}{u} \ll \prod_{p \le x^\varepsilon} \left( 1 + \frac{g_1(p)}{p} \right) \ll \sideset{}{_x}\prod  (g_1) \prod_{x^\varepsilon < p \le x} \left( 1 + \frac{g_1(p)}{p} \right)^{-1}.
\]

From the lower bound hypothesis on $g$ and the construction of the sequence $r$, an integration by parts shows that
\[
\sum_{x^\varepsilon < p \le x} \frac{1}{p} g_1(p) \ge \frac{1}{2}(c-\alpha) \log \frac{1}{\varepsilon} + O(1).
\]
The contribution to $M$ from the terms with $u \le x^\varepsilon$ is
\[
\ll \varepsilon^{(c-\alpha)/2} x(\log x)^{-1} \sideset{}{_x}\prod  (g), \quad x \to \infty.
\]

For the range $x^{1-\varepsilon} < u \le x$, $v \le x^\varepsilon$ and we may invert summations, replacing $(c-\alpha)/2$, as the exponent of $\varepsilon$, by $\alpha/2$.

We are reduced to the estimation of
\[
M_\varepsilon = \sum_{x^\varepsilon < u \le x^{1-\varepsilon}} h_1(u) \sum_{v \le x/u} h_2(v).
\]

Since $h_2$ inherits its properties relative to $g_2$ from $h$, applied to the innersum in $M_\varepsilon$, Satz \ref{wirsing_1967_satz_01_02_02} delivers the asymptotic estimate
\[
\frac{e^{-\gamma \alpha}}{\Gamma(\alpha)} \frac{x/u}{\log(x/u)} \left( \sideset{}{_{x/u}}\prod  (h_2) + o\left( \sideset{}{_{x/u}}\prod  (g_2) \right) \right), \quad x \to \infty,
\]
uniformly for $x^\varepsilon \le u \le x^{1-\varepsilon}$.

Introducing factors $\exp(-p^{-1} h_2(p))$, $\exp(-p^{-1} g_2(p))$, respectively, the ratio $\sideset{}{_{y}}\prod  (h_2) \sideset{}{_{y}}\prod  (g_2)^{-1}$ has an estimate
\[
(B+ o(1)) \exp\left( - \sum_{p \le y} p^{-1} (g_2(p) - h_2(p)) \right), \quad y \to \infty,
\]
with
\[
B = \prod_p \left( \sum_{k=0}^\infty p^{-k} h_2(p^k) \exp(-p^{-1} h_2(p)) \right) \prod_p \left( \sum_{m=0}^\infty p^{-m} g_2(p^m) \right)^{-1} \exp(p^{-1} g_2(p)).
\]

If the series $\sum p^{-1}( g_2(p) - \Re( h_2(p)))$ diverges, then uniformly for $x^\varepsilon \le u \le x^{1-\varepsilon}$,
\[
\sideset{}{_{x/u}}\prod  (h_2) \sideset{}{_{x/u}}\prod (g_2)^{-1} = \sideset{}{_x}\prod  (h_2) \sideset{}{_x}\prod  (g_2)^{-1} + o(1), \quad x \to \infty,
\]
since both product ratios asymptotically vanish.

If the series $\sum p^{-1} (g(p) - \Re(h_2(p)))$ converges, then we may argue as in the proof of Theorem \ref{elliott_2015_multiplicative_function_mean_values_thm_01}.  For each positive real $\tau$, $0 < \tau \le 1$, 
\[
\sum_{x^\tau < p \le x} p^{-1} (g_2(p) - h_2(p)) \to 0, \quad x \to \infty,
\]
and we formally obtain the same asymptotic equality of ratios.

Likewise, there is a representation
\[
(\log y)^{-\alpha} \sideset{}{_{y}}\prod (g_2) = (C + o(1)) \exp\left( \sum_{p \le y} p^{-1} g_2(p) - \alpha \log\log y \right), \quad y \to \infty,
\]
with
\[
C = \prod_p \left( \sum_{m=1}^\infty p^{-m} g_2(p^m) \right) \exp ( - p^{-1} g_2(p)).
\]
An integration by parts shows that for each $\tau$, $0 < \tau < 1$,
\[
\sum_{x^\tau < p \le x} p^{-1} g_2(p) + \alpha \log \tau \to 0, \quad x \to \infty,
\]
so that
\[
(\log(x/u))^{-\alpha} \sideset{}{_{x/u}}\prod  (g_2) = (\log x)^{-\alpha} \sideset{}{_x}\prod  (g_2) + o(1), x \to \infty.
\]

Altogether, the innersum of $M_\varepsilon$ has the estimate
\[
\frac{e^{-\gamma \alpha}}{\Gamma(\alpha)} \frac{x}{u(\log x)^\alpha} \cdot \frac{1}{(\log (x/u))^{1-\alpha}} \left( \sideset{}{_x}\prod  (h_2) + o\left( \sideset{}{_x}\prod  (g_2) \right) \right), \quad x \to \infty,
\]
uniformly for $x^\varepsilon \le u \le x^{1-\varepsilon}$.

The error terms contribute towards $M_\varepsilon$
\[
o\left( \frac{x}{\log x} \sideset{}{_x}\prod  (g) \sum_{x^\varepsilon < u \le x^{1-\varepsilon}} \frac{g_1(u)}{u} \right) = o\left( \frac{x}{\log x} \sideset{}{_x}\prod  (g) \right), \quad x \to \infty,
\]
within which $M_\varepsilon$ has the estimate
\[
\frac{e^{-\gamma}}{\Gamma(\alpha)} \frac{x}{(\log x)^\alpha} \sideset{}{_x}\prod  (h_2) \sum_{x^\varepsilon < u \le x^{1-\varepsilon}} \frac{h_1(u)}{u(\log (x/u))^{1-\alpha}}.
\]

Setting
\[
H_1(y) = \sum_{n \le y} h_1(n)n^{-1}, \qquad G_1(y) = \sum_{n \le y} g_1(n) n^{-1},
\]
an integration by parts gives a representation
\[
\sum_{x^\varepsilon < u \le x^{1-\varepsilon}} \frac{h_1(u)}{u(\log(x/u))^{1-\alpha}} = \frac{H_1(x^{1-\varepsilon})}{(\varepsilon \log x)^{1-\alpha}} - \frac{H_1(x^\varepsilon)}{((1-\varepsilon)\log x)^{1-\alpha}} - (1-\alpha)\int_{x^\varepsilon}^{x^{1-\varepsilon}} \frac{H_1(y)}{y(\log(x/u))^{2-\alpha}} \, dy,
\]
provided $x^\varepsilon$, $x^{1-\varepsilon}$ are not positive integers, a situation that we may avoid by choosing a slightly larger value of $x$.

According to Theorem \ref{elliott_2015_multiplicative_function_mean_values_thm_01},
\[
H_1(y) = \left( \sideset{}{_{y}}\prod  (h_1) \sideset{}{_{y}}\prod  (g_1)^{-1} + o(1) \right) G_1(y), \quad y \to \infty,
\]
where, as above, we may replace the products $\sideset{}{_{y}}\prod $ by $\sideset{}{_x}\prod $, uniformly for $x^\varepsilon \le y \le x^{1-\varepsilon}$, $x \to \infty$.

As a consequence,
\[
\sum_{x^\varepsilon < u \le x^{1-\varepsilon}} \frac{h_1(u)}{u(\log (x/u))^{1-\alpha}} =  \left( \sideset{}{_{x}}\prod  (h_1) \sideset{}{_{x}}\prod  (g_1)^{-1} + o(1) \right) \sum_{x^\varepsilon < u \le x^{1-\varepsilon}} \frac{g_1(u)}{u(\log(x/u))^{1-\alpha}}, \quad x \to \infty.
\]
Once again, the argument is expedited by considering $2G_1(x) + \Re (H_1(x))$, $2G_1(x) + \Im( H_1(x))$.

Rewinding,
\begin{align*}
M_\varepsilon & = \sideset{}{_x}\prod (h_1) \sideset{}{_x}\prod (g_1)^{-1} \frac{e^{-\gamma\alpha}}{\Gamma(\alpha)} \frac{x}{(\log x)^\alpha} \sideset{}{_x}\prod (h_2) \sum_{x^\varepsilon < u \le x^{1-\varepsilon}} \frac{g_1(u)}{u(\log(x/u))^{1-\alpha}} + o\left( \frac{x}{\log x} \sideset{}{_x}\prod (g) \right) \\
& = \sideset{}{_x}\prod (h) \sideset{}{_x}\prod  (g)^{-1} \sum_{x^\varepsilon < u \le x^{1-\varepsilon}} \frac{e^{-\gamma\alpha}}{\Gamma(\alpha)} \frac{x}{\log(x/u)} \sideset{}{_{x/u}}\prod (g_2) \frac{g_1(u)}{u} + o\left( \frac{x}{\log x} \sideset{}{_x}\prod (g) \right) \\
& = \sideset{}{_x}\prod (h) \sideset{}{_x}\prod (g)^{-1} \sum\limits_{x^\varepsilon < u \le x^{1-\varepsilon}} g_1(u) \sum_{v \le x/u} g_2(v) + o\left( \frac{x}{\log x} \sideset{}{_x}\prod (g) \right) \\
& = \sideset{}{_x}\prod (h) \sideset{}{_x}\prod  (g)^{-1} \sum_{n \le x} g(n) + O\left( \varepsilon^\nu \sum_{n \le x} g(n) \right)
\end{align*}
with $\nu = \min((c-\alpha)/2, \alpha/2)$ and, for all sufficiently large values of $x$, an implied constant independent of $\varepsilon$.

A similar estimate holds for $M$.

Letting $x \to \infty$, $\varepsilon \to 0+$ completes the proof.

%@@@@@@@@@@@@@@@@@@@@@@@@@@@@@@@@@@@@@@@@@@@@@@@@@@@@@@@@@

\section{Proof of Theorem \ref{elliott_2015_multiplicative_function_mean_values_thm_04}}\label{elliott_2015_multiplicative_function_mean_values_sec_06}

\noindent \emph{Case (i)}.  From the assumption that the series $\sum p^{-1} (g(p) - \Re (h(p)p^{it}))$ converges, for each positive $\delta$ the series taken over the primes $p$ for which $g(p) - \Re (h(p)p^{it})> \delta$ also converges.

On the remaining primes
\[
| g(p) - h(p)p^{it}|^2 \le 2g(p) (g(p) - \Re (h(p)p^{it})) \le 2\beta\delta.
\]
The values of $h(p)p^{it}$ lie in a box about the real axis, with corners at $(-(2\beta\delta)^{1/2}, \pm (2\beta\delta)^{1/2})$, $(\beta + (2\beta\delta)^{1/2}, \pm (2\beta\delta)^{1/2})$, and area $2(2\beta\delta)^{1/2}(\beta + 2(2\beta\delta)^{1/2})$.

Assuming that $\delta$ is sufficiently small and, in particular, that $2(2\beta\delta)^{1/2} \le \beta$, this is a region of the type $\Delta(\tau)$ with an average radius
\[
\frac{1}{2\pi} \int_0^{2\pi} w(\theta) \, d\theta \le \left( \frac{1}{2\pi} \int_0^{2\pi} w(\theta)^2 \, d\theta \right)^{1/2} \le \left( 4\pi^{-1}(2\beta^3\delta)^{1/2} \right)^{1/2}
\]
that can be fixed at a value as small as desired.

We may follow the proof of Theorem \ref{elliott_2015_multiplicative_function_mean_values_thm_03}, first selecting a subsequence of primes $r$ for which $(\log x)^{-1} \sum_{r \le x} r^{-1}g(r)\log r \to \alpha$, $x \to \infty$, then removing from that subsequence those primes for which $h(p)p^{it}$ does not belong to a region $\Delta(\alpha)$ defined by a value of $\delta$ that satisfies $4\pi^{-1}(2\beta^3\delta)^{1/2} < \alpha^2$.

The removal of these exceptional primes does not affect the existence or the value of the asymptotic limit for $(\log x)^{-1} \sum_{r \le x} r^{-1} g(r) \log r$.

The upshot is an asymptotic estimate
\[
\sum_{n \le x} h(n)n^{it} = \prod_{p \le x} \left( 1+ h(p)p^{it-1} + \cdots \right) \left( \sideset{}{_x}\prod  (g) \right)^{-1} \sum_{n \le x} g(n) + o\left( \sum_{n \le x} g(n) \right), \quad x \to \infty.
\]

We would like to integrate by parts and remove the weight $n^{it}$ from $h(n)n^{it}$, but have insufficient control over the values of the function $h$.  Since, in some sense, we are considering the ratio $h(n)n^{it}(g(n))^{-1}$, at an appropriate moment we switch the weight $n^{it}$ from $h$ to $g$ and consider the ratio $h(n)(g(n)n^{-it})^{-1}$.

Following the argument for Theorem \ref{elliott_2015_multiplicative_function_mean_values_thm_03}, the study of the sum $\sum_{n \le x} h(n)$ is reduced to that of
\[
\widetilde{M}_\varepsilon = \sum_{x^\varepsilon < u \le x^{1-\varepsilon}} h_1(u) \sum_{v \le x/u} h_2(v),
\]
where Theorem \ref{elliott_2015_multiplicative_function_mean_values_thm_03} is applicable to the pair $h_2(n)n^{it}$, $g_2(n)$.  There is a corresponding estimate
\[
\sum_{n \le y} h_2(n)n^{it} = L(\log y) \sum_{n \le y} g_2(n) + o\left( \sum_{n \le y} g_2(n) \right), \quad y \to \infty,
\]
with
\[
L(\log y) = \prod_{p \le y} \left( 1 + h_2(p)p^{it-1} + \cdots \right) \left( \sideset{}{_{y}}\prod  (g_2) \right)^{-1}, \quad y \ge 2.
\]

Set
\[
H_2(y) = \sum_{n \le y} h_2(n)n^{it}, \qquad G_2(y) = \sum_{n \le y} g_2(n), \quad y \ge 1/2.
\]
An integration by parts gives a representation
\[
\sum_{n \le y} h_2(n) = y^{-it} H_2(y) + it \int_{1/2}^y w^{-it-1}H_2(w) \, dw,
\]
provided $y$ is not an integer.  Since $G_2(w) \ll w(\log w)^{-1} \prod_w (g_2)$, $w \ge 2$,
\begin{align*}
\int_2^x w^{-1}G_2(w) \, dw & \ll \sideset{}{_x}\prod  (g_2) \int_2^x (\log w)^{-1} \, dw \\
& \ll x(\log x)^{-1} \sideset{}{_x}\prod  (g_2) \ll G_2(x), \quad x \ge 2.
\end{align*}
Hence
\[
\sum_{n \le y} h_2(n) = y^{-it} L(\log y)G_2(y) + it\int_2^y w^{-it-1} L(\log w)G_2(w) \, dw + o(G_2(y)), \quad y \to \infty.
\]

As in the proof of Theorem \ref{elliott_2015_multiplicative_function_mean_values_thm_03}, within an acceptable error $L(\log w)$, for $y^\varepsilon \le w \le y$, may be replaced by $L(\log y)$ and factored out of the representation:
\[
\sum_{n \le y}h_2(n) = L(\log y)\left( y^{-it}G_2(y) + it\int_2^y w^{-it-1} G_2(w) \, dw \right) + o(G_2(y)), \quad y \to \infty.
\]

We appeal to the asymptotic estimate
\[
G_2(x) = (1 + o(1)) \frac{e^{-\gamma \alpha}}{\Gamma(\alpha)} \frac{x}{\log x} \sideset{}{_x}\prod  (g_2), \quad x \to \infty,
\]
vouchsafed by Satz \ref{wirsing_1967_satz_01_01}.  Once again, as for Theorem \ref{elliott_2015_multiplicative_function_mean_values_thm_03}, we employ the slow oscillation of the function $\sideset{}{_x}\prod  (g_2) (\log x)^{-\alpha}$ to obtain a representation
\begin{align*}
\sum_{n \le y} h_2(n) & = \frac{e^{-\gamma\alpha}}{\Gamma(\alpha)} L(\log x) \frac{\sideset{}{_x}\prod  (g_2)}{(\log x)^\alpha} \left( \frac{y^{1-it}}{(\log y)^{1-\alpha}} + it \int_2^y \frac{w^{-it}}{(\log w)^{1-\alpha}} \, dw \right) + o(G_2(y)) \\
& = \frac{e^{-\gamma\alpha}}{\Gamma(\alpha)} L(\log x) \frac{\sideset{}{_x}\prod  (g_2)}{(\log x)^\alpha} \frac{y^{1-it}}{(1-it)(\log y)^{1-\alpha}} + o(G_2(y)),
\end{align*}
uniformly for $x^\varepsilon \le y \le x$, as $x \to \infty$; stepping from $w$ to $y$ to $x$.

Accordingly,
\[
\widetilde{M}_\varepsilon = \frac{e^{-\gamma\alpha}}{\Gamma(\alpha)} \frac{x^{1-it}}{1-it} L(\log x) \frac{\sideset{}{_x}\prod  (g_2)}{(\log x)^\alpha} \sum_{ x^\varepsilon < u \le x^{1-\varepsilon}} \frac{h_1(u)u^{it}}{u(\log x/u)^{1-\alpha}} + o(G(x)), \quad x \to \infty.
\]

We may now formally follow the argument for Theorem \ref{elliott_2015_multiplicative_function_mean_values_thm_03}, the r\^ole of $h_1(n)$ there here played by $h_1(n)n^{it}$, although on a slightly different set of primes.  Eventually only the extra factor $x^{-it}(1-it)^{-1}$ remains.

\noindent \emph{Case (ii)}.  The series $\sum p^{-1}( g(p) - \Re (h(p)p^{it}))$ diverges for every real $t$.  The partial sums of this series are non-decreasing in $x$ and continuous in $t$.  Divergence of the series is uniform on every compact interval $|t| \le T$ and Theorem \ref{elliott_2015_multiplicative_function_mean_values_thm_03} follows from an application of Theorem \ref{elliott_2015_multiplicative_function_mean_values_thm_06}.

\noindent \emph{Remark}.  Under the hypothesis of Case (i) the series $\sum p^{-1} |g(p) - h(p) p^{it}|^2$ converges.  The series $\sum p^{-1} |g(p) - |h(p)||^2$ and $\sum p^{-1} g(p) | 1 - e^{i\theta_p} p^{it} |^2$, where $h(p) = |h(p)| e^{i\theta_p}$, then also converge.

Suppose further that, for some positive integer $k$, $h(p)^k$ is real.  The inequality $|1-z^k| \le k|1-z|$, valid for every $z$ in the complex unit disc, guarantees the series $\sum p^{-1} g(p) |1-p^{2ikt}|^2$ to converge.

In the present circumstances $\sum_{p \le x} p^{-1} g(p) \ge (c+o(1)) \log\log x$ as $x \to \infty$ and an application of Lemma 15 from Elliott and Kish \cite{elliottkish2013harmonicbasic} shows that $t=0$.

A simple example is given by $h(n) = g(n)\chi(n)$, where $\chi$ is a Dirichlet character.

The argument of this remark may be given a topological aspect by defining a metric $\sigma(f,g) = ( \sum p^{-1} |f(p) - g(p)|^2)^{1/2}$ on equivalence classes of multiplicative functions that coincide of the primes, and restricting study to those functions $g$ whose distance $\sigma(g, g_0)$ to a fixed multiplicative function $g_0$ is defined, i.e. finite.  The topological space of complex-valued multiplicative functions is in this manner locally metrised and correspondingly disconnected.

%@@@@@@@@@@@@@@@@@@@@@@@@@@@@@@@@@@@@@@@@@@@@@@@@@@@@@@@@@@

\section{Proof of Theorem \ref{elliott_2015_multiplicative_function_mean_values_thm_06}}\label{elliott_2015_multiplicative_function_mean_values_sec_07}

We assume the new, weaker restraints upon $g$.  If $g$ is exponentially multiplicative, i.e. $g(p^k) = g(p)^k/k!$, and $|g(p)| \le \beta$, then for any $\gamma$ the series $\sum_{p, k \ge 2} p^{-k} |g(p^k)| (\log p^k)^\gamma$ converges, so that Theorem \ref{elliott_2015_multiplicative_function_mean_values_thm_04} is applicable.  Indeed, for such functions the original exposition of Elliott and Kish, \cite{elliottkish2013harmonicmaass} Theorem 2, already contains a proof.

In general, we define an exponentially multiplicative function $g_1$ by $g_1(p) = g(p)$, and a complementary multiplicative function $g_2$ by Dirichlet convolution:  $g = g_1* g_2$.

Calculation with Euler products shows that $g_2(p)=0$ and for $k \ge 2$,
\[
g_2(p^k) = \sum_{r=0}^k (r!)^{-1} (-g(p))^r g(p^{k-r}).
\]
In particular,
\[
|g_2(p^k)| \le \sum_{r=0}^k (r!)^{-1} \beta^r |g(p^{k-r})|, \quad k \ge 2.
\]
As a consequence
\[
\sum_{p, k \ge 2} p^{-k} |g_2(p^k)| \le \sum_{r=0}^\infty (r!)^{-1} \beta^r \sum_{p,k \ge 2} p^{-k} |g(p^{k-r})|
\]
\[
\le \left( \tfrac{3}{2} \beta^2 + \tfrac{1}{4}\beta^3 \right) \sum p^{-2} + \left( 1 + \tfrac{1}{2} \beta^2 \right) \sum_{p, k \ge 2} p^{-k} |g(p^k)|,
\]
and converges.

Moreover,
\[
\sum_{p^k \le y} |g_2(p^k)| \le \sum_{r=0}^\infty (r!)^{-1} \beta^r \sum_{p^k \le y, k \ge 2} |g(p^{k-r})|
\]
\[
\ll \sum_{r=0}^\infty (r!)^{-1} \beta^r y (\log y)^{-1} \ll y(\log y)^{-1}
\]
uniformly for $y \ge 2$.

We may apply Lemma \ref{elliott_2015_multiplicative_function_mean_values_lem_01} and obtain for $|g_2|$ the uniform estimate
\[
\sum_{n \le y} |g_2(n)| \ll y(\log y)^{-1}, \quad y \ge 2.
\]

With $\delta$ a real number to be chosen presently in the range $0 < \delta < 1$,
\[
\rho = \exp\left( -\frac{c}{c+\beta} \lambda \right) + T^{-1/2},
\]
as in the statement of Theorem \ref{elliott_2015_multiplicative_function_mean_values_thm_05}, we define $w = \exp(\rho^\delta \log x)$, so that $w$ is effectively a function of $x$ for $x \ge 2$.

It is convenient to note that we may assume $\rho^\delta \le 1/2$, otherwise Theorem \ref{elliott_2015_multiplicative_function_mean_values_thm_06} follows directly from Lemma \ref{elliott_2015_multiplicative_function_mean_values_lem_01}.

Moreover, provided $2\delta\beta c < c+ \beta$ and $Y$ does not exceed a certain fixed power of $x$, which we may likewise assume, $Y \le w$.  For otherwise
\[
\log x/\log Y \le \rho^{-\delta} \le \exp\left( \frac{\delta c}{c + \beta} \lambda \right)
\]
\[
\ll \exp\left( \frac{\delta c }{c + \beta} \sum_{Y < p \le x} 2p^{-1} |g(p)| \right) \ll (\log x/\log Y)^{2\delta\beta c/(c+\beta)}.
\]
In particular, uniformly for $w < y \le x$,
\[
\min\limits_{|t| \le T} \sum_{Y<p \le y} p^{-1} ( |g(p)| - \Re( g(p)p^{it})) \ge \lambda - 2\sum_{w < p \le x} p^{-1} |g(p)| \ge \lambda + 2\delta\beta \log \rho + O(1).
\]

Applied to $g_1$ over the same range of $y$-values, Theorem \ref{elliott_2015_multiplicative_function_mean_values_thm_05} delivers an estimate
\begin{align*}
\sum_{n \le y} g_1(n) & \ll \frac{y}{\log y} \sideset{}{_{y}}\prod  (|g_1|) \left( \exp\left( -\frac{c\lambda}{c + \beta} \right) \rho^{-2\delta\beta c/(c+\beta)} + T^{-1/2} \right) \\
& \ll \frac{y}{\log y} \sideset{}{_{y}}\prod  (|g_1|) \rho^{1-2\delta\beta c/(c+\beta)},
\end{align*}
this last step somewhat wasteful.

We decompose the mean-value of $g$ into two sums:
\[
\sum_{n \le x} g(n) = \sum_{b \le x/w} g_2(b) \sum_{a \le x/b} g_1(a) + \sum_{a < w} g_1(a) \sum_{x/w < b \le x/a} g_2(b).
\]

The first doublesum is
\[
\ll \sum_{b \le x/w} |g_2(b)| x b^{-1} (\log(x/b))^{-1} \sideset{}{_{x/b}}\prod  (|g_1|) \rho^{1-2\delta\beta c/(c+\beta)}
\]
\[
\ll x(\log x)^{-1} \sideset{}{_x}\prod  (|g|) \rho^{1-2\delta\beta c/(c+\beta) - \delta}.
\]

The second doublesum is
\[
\ll \sum_{a < w} |g_1(a)| xa^{-1} (\log (x/a))^{-1}
\]
and $w \le x^{1/2}$, so that the bound does not exceed a constant multiple of
\[
x(\log x)^{-1} \prod_{p \le w} ( 1 + p^{-1}|g(p)|) \ll x(\log x)^{-1} \sideset{}{_x}\prod  (|g|) \exp\left( -\sum_{w < p \le x} p^{-1} |g(p)| \right).
\]
According to the lower bound hypothesis on $|g(p)|$ in Theorem \ref{elliott_2015_multiplicative_function_mean_values_thm_05}, still in force in Theorem \ref{elliott_2015_multiplicative_function_mean_values_thm_06}, noting that $w \ge Y$,
\[
\sum_{w < p \le x} p^{-1} |g(p)| \ge c\sum_{w < p \le x} p^{-1} + O(1) \ge -\delta c \log \rho + O(1).
\]

Altogether,
\[
\sum_{n \le x} g(n) \ll \frac{x}{\log x} \sideset{}{_x}\prod  (|g|) ( \rho^{1-\delta c_0} + \rho^{\delta c} )
\]
with $c_0 = 2\beta c(c+\beta)^{-1}+1$.

We choose $\delta$ to satisfy $1-\delta c_0 = \delta c$.  The earlier condition $2\delta \beta c < \beta + c$ is amply satisfied, $c_0$ increases with $\beta$ and $\delta c$ descends to a limiting value $c(3c+1)^{-1}$.

%@@@@@@@@@@@@@@@@@@@@@@@@@@@@@@@@@@@@@@@@@@@@@@@@@@@@@@@@@@

\section{Concluding Remarks}\label{elliott_2015_multiplicative_function_mean_values_sec_08}

The hypothesis on $|g|$ in Theorem \ref{elliott_2015_multiplicative_function_mean_values_thm_06} remains essentially weaker than that on $g$ in Theorem \ref{elliott_2015_multiplicative_function_mean_values_thm_04}.  What might a best-possible condition on $g$ be in order to guarantee the validity of Theorem \ref{elliott_2015_multiplicative_function_mean_values_thm_04}?

Likewise, what might the weakest hypothesis on $g$ be in order to guarantee the validity of the lower bound in Theorem \ref{elliott_2015_multiplicative_function_mean_values_thm_02}?

\bibliographystyle{amsplain}
\bibliography{MathBib}

\end{document}